\let\documentclass\relax\documentclass % Для тупого арХива
\input AHTOH-E.STY
\def\r{{\bf\=r}}
\def\NO{{\bf NO}}
\def\rk{\={\rm rk}}    %переопределён
\def\rank {{\rm rank}} %переопределён слегка

\UDC{
512.543.12%  %Свободные группы
+512.543.52% %Свободные произведения
+512.543.14% %Группы с конечным числом образующих
+512.544.23% %Структурные свойства групп
}


\title{%
Intersections of subgroups in virtually free groups
and virtually free products
}
\author{%
Anton A. Klyachko
\quad
Anastasia N. Ponfilenko
}
\address{
\myAddressW
\quad
stponfilenko@gmail.com
}

\grantsFirst{\RFBR19-01-00591}

\abstract{%
%\hfuzz12pt
%Our results implies, in particular, 
This note contains a (short) proof of
the following generalisation of
the Friedman--Mineyev theorem (earlier known as the
Hanna Neumann conjecture): if $A$ and $B$ are nontrivial free subgroups
of a virtually free group containing a free subgroup of index $n$, then
%\newline
%\centerline
{$\rank(A\cap B)-1\le n\cdot(\rank(A)-1)\cdot(\rank(B)-1)$.}
In addition, we obtain a virtually-free-product analogue of this 
result.}

%%%%%%%%%%%%%%%%%%%%%%%%%
\s 1.
Introduction

The Hanna Neumann Conjecture (1957) proven by
Mineyev ([Mi12a], [Mi12b]) and Friedman [Fr14] asserts that
\disp{\sl
for any nontrivial subgroups $A$ and $B$ of a free group,
%\newline
%\centerline{
$\rank(A\cap B)-1\le (\rank(A)-1)\cdot(\rank(B)-1)$.}
%}

We obtain the following generalisation.

\Th 1.
For any nontrivial free subgroups $A$, $B$,
and $F$ of any group $G$,
$$
\rank(A\cap B)-1\le|G{:}F|\cdot(\rank(A)-1)\cdot(\rank(B)-1).
\eqno{(1)}
$$

This inequality can be understood in the sense of cardinal arithmetics,
but it is nontrivial only when all
three value on the right-hand side --- the rank of $A$, rank
of~$B$, and index of~$F$ --- are finite.

It has previously been known that
$$
\eqalign{
\rank(A\cap B)-1&\le6|G{:}F|(\rank(A)-1)(\rank(B)-1)
\qquad\qquad\qquad\qquad\hbox{[Za14]};
%\eqalignno{\hbox{[Za14]};}
\cr
\rank(A\cap B)-1&\le|G{:}F|^2(\rank(A)-1)(\rank(B)-1)+|G{:}F|-1
\qquad\hbox{[ASS15]}.
}
$$
Of course, the estimate from [ASS15] is asymptotically worse than
the estimate from [Za14], but it is better for some small values
of the index.
Theorem 1 improves both these inequalities; and
no further improvement
is possible:
\disp{\sl
\narrower \narrower
for any $k,l,n\in\N$, there
exists a group $G$ containing free subgroups $A$, $B$, and $F$ such that
%\newline
$\rank(A)=k$, \ $\rank(B)=l$, \ $|G{:}F|=n$, and
inequality {\rm(1)} is an equality.
}%
Indeed,
consider an epimorphism $\phi\:x\mapsto\bigl(\alpha(x),\beta(x)\bigr)$
from the free group $F$ of rank two
onto the free abelian group~$\Z\oplus\Z$
and take the groups
$$
A=\alpha^{-1}\bigl((k-1)\Z\bigr),
\quad
B_0=\beta^{-1}\bigl((l-1)\Z\bigr),
\quad
G=F\times(\Z/n\Z)\supset
B=\Bigl\{\Bigl(b,\;\beta(b)/(l-1)\Bigr)\;\Bigm|\;b\in B_0\Bigr\}.
$$
Thus,
$$
B\cap F=\beta^{-1}\bigl(n(l-1)\Z\bigr)
\qqbox{and}
B\cap A=\phi^{-1}\bigl((k-1)\Z\times n(l-1)\Z\bigr).
$$
Clearly, $|G{:}F|=n$. Moreover, $|F{:}A|=k-1$, $|F{:}B_0|=l-1$ and
$|F{:}B\cap A|=n(l-1)(k-1)$. The ranks of these subgroups are
$k$, $l$, and $n(l-1)(k-1)+1$, respectively,
by the Schreier formula:
$\rank(H)-1=|F{:}H|(\rank(F)-1)$ (which is valid for any subgroup
$H$ of finite index in a free group $F$).
It remains to note that
$\rank(B)=\rank(B_0)$
because the projection $\Bigl(b,\;\beta(b)/(l-1)\Bigr)\mapsto b$
is an isomorphism from $B$ to $B_0$.

\medskip
(We note parenthetically that Theorem 1.8 stated in [Mi12b]
without proof would imply estimate (1) without the factor $|G{:}F|$;
obviously that theorem contains some misprints.)

The following result can be considered as a generalisation of Theorem 1.
Recall that a \emph{left-orderable group} is a group admitting
a linear (=~total) order such that $x\le y\imp zx\le zy$ for any
elements $x,y,z$.

\Th 2.
Suppose that a group $G$ has a finite-index subgroup
$F=\zvezda_{i\in I}G_i$
which is a free product of left-orderable groups $G_i$.
Then, for
any nontrivial free subgroups~$A$~and~$B$ of $G$ trivially
intersecting all subgroups conjugate to $G_i$,
the following inequality holds
$$
\rank(A\cap B)-1\le|G{:}F|(\rank(A)-1)(\rank(B)-1).
$$

For $F=G$ this assertion was proven in [AMS14] (see also [Iv17]).

The authors thank A.~O.~Zakharov for reading a previous version of this 
text and suggesting corrections and refinements. 

%%%%%%%%%%%%%%%%%%%%%%%%%
\s 2.
Tools

Throughout this paper, the word \emph{graph} means a directed
graph; loops and multiple edges are allowed. A \emph{path} in a graph and
the
\emph{connectedness} of a graph are defined naturally (ignoring
directions). The \emph{reduced rank} $\r(D)$ of a finite graph $D$
is defined as $\r(D)\:=\sum\limits_K\max(0,-\chi(K))$, where
the sum is over all connected components $K$ of $D$, and $\chi(K)$
is the \emph{Euler characteristic} of the graph~$K$, i.e. the difference
between the number of vertices and edges.

\noindent
We call a set $E$ of edges of a graph $D$
\emph{maximal essential} if
$\r(D\setminus E)=\r(D)-|E|=0$.
In other words, a set $E$ of edges of a graph $D$
is maximal essential if $D\setminus E$ is inclusion-maximal
subgraph of $D$ whose each component is homotopic to either a point or
a circle.

We say that a graph
is \emph{ordered} if the set of its of edges is partially
ordered such that the order on
each connected component is linear.
An edge $e$ of an ordered forest is called
\emph{order-essential} if it lies in a
bi-infinite simple path consisting of
edges not exceeding $e$.

\noindent
An action of a group on a graph is called
\-
\emph{cocompact}
if the number of orbits of vertices and edges are finite;
\-
\emph{free}
if the stabiliser of each vertex is trivial
(and, therefore, the stabilisers of edges are also trivial);
\-
\emph{free on edges}
if the stabiliser of each edge is trivial.

\proclaim Mineyev's essential-edge theorem
\rm([Mi12b], Theorem 1.6).
%теорема 5 в препринте.
Suppose that a group~$G$ acts
on an ordered
forest $T$
freely, cocompactly,
and order-preservingly.
Then the set of orbits of order-essential
edges is a maximal essential set in
the quotient graph~$T/G$. In particular, the reduced rank $\r(T/G)$
of the quotient graph
equals the number of
orbits of order-essential edges.

\proclaim{Free-rank lemma}.
If a
free finitely generated group $A$ acts
freely cocompactly and order-preservingly on an ordered forest
$L$
consisting of $n$ trees,
then the number of orbits of order-essential edges equals
$n\cdot\rk(A)$.
\newline
\rm Henceforth, $\rk(A)\:=\max(\rank(A)-1,0)$ is the \emph{reduced
rank} of a free group $A$.

\Proof
Let $\NO(G,\Gamma)$ denote the number of $G$-orbits of
order-essential
edges in an ordered graph $\Gamma$
(on which a group~$G$ acts preserving the order).

\Case 1: $n=1$.
In this case, (as
was noted in [Mi12b], Lemma 1.1)
the assertion
follows immediately from Mineyev's essential-edge theorem because
$\r(T/A)=\rk(A)$ if $T$ is a tree.

\Case 2: the action of $A$ on the set of
components of $L$ is transitive.
Let $T$ be a tree (component) of the forest~$L$.
Then
$
\NO(A,L)
\buildrel t\over=
\NO(\St(T),T)
\buildrel1\over=
\rk(\St(T))
\buildrel S\over=
|A:\St(T)|\cdot\rk(A)
=
n\cdot\rk(A),
$
where equality $\buildrel t\over=$ follows from the transitivity
of the action on the
set of components; equality $\buildrel 1\over=$ is
Case 1; equality~$\buildrel S\over=$ is the
Schreier formula, and the last equality is the orbit-stabiliser
theorem.

\Case 3: general case.
Suppose that $L=P_1\sqcup\dots\sqcup P_k$ and,
on each (invariant) forest $P_i$ consisting of $l_i$
trees, the action
is transitive. Then
$$
\NO(A,L)=\NO(A,P_1)+\dots+\NO(A,P_k)
\buildrel 2\over=
l_1\cdot\rk(A)+\dots+l_k\cdot\rk(A)
=
(l_1+\dots+l_k)\cdot\rk(A)
=
n\cdot\rk(A),
$$
where equality $\buildrel2\over=$ is Case 2.
This completes the proof.

\goodbreak
\medskip

The following simple lemma can be found, e.g., in [Za14] (Lemma 2).

\proclaim Orbit-intersection lemma.
Suppose that $A$ and $B$ are subgroups of a group $G$
acting freely on
a set~$X$ containing an
$A$-invariant
subset $Y\subseteq X$ and a $B$-invariant
subset $Z\subseteq X$.
Then
$$
\hbox{\(the number of $(A\cap B)$-orbits in $Y\cap Z$\)}
\leqslant
\hbox{\(the number of $A$-orbits in $Y$\)}
\cdot
\hbox{\(the number of $B$-orbits in $Z$\)}.
$$

\proclaim Induced-action lemma.
Suppose that a group $G$ has a
subgroup $F$ of a finite index $n$,
and $F$
acts on an ordered tree $T$ preserving the order.
Then $G$
can
order-preservingly
act on an
ordered forest consisting of $n$ trees such that
the stabilisers of vertices and edges are
conjugate to the stabilisers of vertices and edges
under the initial action of $F$
on $T$.

\Proof
Let $S\ni1$ be a system of representatives of the left cosets of $F$ in $G$
(i.e.  $|S|=n$). Thus, each element $g\in G$ decomposes uniquely into a
product $g={\bf s}(g){\bf f}(g)$ of an element~${\bf s}(g)\in S$ and an
element ${\bf f}(g)\in F$.

Take the ordered forest $L=\bigcup\limits_{s\in S} sT$ consisting of
$n$ copies $sT$ of the ordered tree $T$ (edges from different
copies are incomparable) and consider the usual induced
action of $G$ on $L$:
\quad $g\o st\:={\bf s}(gs)\Bigl({\bf f}(gs)\o t\Bigr)$.
Clearly, this action satisfies all requirements.

\proclaim Invariant-forest lemma.
If a finitely generated group $G$ acts on a
forest~$L$ with finitely many connected components, then any finite
set $X$ of vertices of $L$ is contained in a $G$-invariant
subforest~$L_X\supseteq X$ such that its intersection with each
component of $L$ is connected, and the action of~$G$ on $L_X$ is cocompact.

\Proof
For each component $T$ of $L$,
we choose a finite set $S$ of generators of the stabiliser
of~$T$ (the stabiliser is finitely generated as it is a
finite-index subgroup of the
finitely generated group~$G$). Now, we
\-
join all points of $X\cap T$ by (shortest) paths;
\-
join the obtained tree $R$ by paths with the trees $s^{\pm1}R$ for all
$s\in S$ and add these paths to $R$.
\enditem
Finally, we
add to the obtained finite forest $R'\supseteq X$ all its shifts $gR'$,
where $g\in G$. Clearly, we obtain a $G$-invariant 
forest~$R''=\bigcup\limits_{g\in G}gR'$, and the action of $G$ on $R''$ is
cocompact.

Let us verify that the intersection $R''\cap T$ is connected for
each component $T$ of $L$. Indeed, the tree
$R'\cap T$ is joined by paths
with trees $s^{\pm1}R'\cap T$,
hence, the tree $gR'\cap T$ is joined with the trees
$gs^{\pm1}R'\cap T$ for all $g\in\St(T)$ and $s\in S$;
in particular, the trees $gR'\cap T$
and $g'R'\cap T$
lie in the same component of the forest $R''\cap T$,
where the length of the element $g'=gs^{\pm1}\in\St(T)$ (with respect to
the generating set $S$) is less than the length of $g$. An obvious
induction completes the proof.

%%%%%%%%%%%%%%%%%%%%%%%%%
\s 2.
Proof of Theorems 1 and 2

Suppose that a group $G$ has a subgroup $F$
which is either free or at least
decomposes
into a free product of left-orderable groups.
As is known, a free group acts freely on some
tree $T$, and a free product
$F=\zvezda_{i\in I}G_i$
acts on some tree $T$ in such a way that
 the stabiliser of each vertex is conjugate to one of the
factors $G_i$.

The tree $T$ can be ordered:
the order on the set of edges of $T$ is induced by a
left-invariant order on the group $F$ (which is known
to exist [Vi49], [D\v S14]). Thus, the action of $F$ on
$T$ preserves the order and is free on edges
(in both cases).

By the induced-action lemma, the group $G$
acts on an ordered forest $L$ freely on edges.
Moreover, the action of groups $A$ and $B$ on $L$ are free.
By the invariant-forest lemma, we choose an $A$-invariant
subforest ${\cal A}\subseteq L$ and a $B$-invariant
subforest~${\cal B}\subseteq L$
such that
\-
the actions of $A$ and $B$ on
$\cal A$ and $\cal B$ are cocompact
(and, hence, the action of $A\cap B$ on $\cal A\cap\cal B$
is cocompact too by the orbit-intersection lemma);
\-
the intersection of 
each component of $L$ with
each of the
forests $\cal A$, $\cal B$, and $\cal A\cap\cal B$
is nonempty and connected.

\enditem
Setting in the orbit-intersection lemma
$$
\eqalign{
X=\{\hbox{edge of $L$}\},
\qquad
Y&=\{\hbox{order-essential edge of $\cal A$}\},
\cr
Z&=\{\hbox{order-essential edge of $\cal B$}\}.
}
$$
and noting that each order-essential edge
of $\cal A\cap B$ is order-essential in
both $\cal A$ and $\cal B$,
we obtain
$$
\NO(A\cap B,{\cal A\cap B})\le \NO(A,{\cal A})\cdot\NO(B,{\cal B}).
$$ 
By the free-rank lemma, the left-hand side of this inequality equals
$n\cdot\rk(A\cap B)$, and the right-hand side equals
$n^2\cdot\rk(A)\cdot\rk(B)$.
Cancelling out $n$, we obtain (1).

%%%%%%%%%%%%%%%%%%%%%%%%%
\References

[Vi49]
A. A. Vinogradov,
On the free product of ordered groups,
Mat. Sb. (N.S.), 25(67):1 (1949), 163-168.

%[LM17]
%V. Lebed, A. Mortier,
%A translation of A.A.Vinogradov's "On the free product of ordered groups",
%arXiv:1703.05781.

%[Ши47]
%Е. П. Шимбирева,
%К теории частично упорядоченных групп,
%Матем. сб., 20(62):1 (1947), 145-178.

[AMS14]
Y. Antolin, A. Martino, and I. Schwabrow,
Kurosh rank of intersections of subgroups of free products of
right-orderable groups,
Mathematical Research Letters, 21:4 (2014), 649-661.
\arXiv 1109.0233

[ASS15]
V. Ara\'ujo, P. V. Silva, M. Sykiotis,
Finiteness results for subgroups of finite extensions,
J. Algebra, 423 (2015), 592-614.
\arXiv 1402.0401

[D\v S14]
W. Dicks, Z. \v Suni\'c,
Orders on trees and free products of left-ordered groups,
arXiv:1405.1676.

[Fr14]
J. Friedman,
Sheaves on graphs, their homological invariants, and a proof of the
Hanna Neumann conjecture. With an appendix by Warren Dicks,
Mem. Amer. Math. Soc. 233:1100 (2014).
\arXiv 1105.0129

[Iv17]
S. V. Ivanov,
Intersecting free subgroups in free products of left ordered groups,
Journal of Group Theory, 20:4 (2017), 807-821.
\arXiv 1607.03010

[Mi12a]
I. Mineyev,
Submultiplicativity and the Hanna Neumann conjecture,
Ann. Math., 175 (2012), 393-414.

[Mi12b]
I. Mineyev,
Groups, graphs, and the Hanna Neumann conjecture,
J. Topol. Anal., 4:1 (2012), 1-12.

[Za14]
A. Zakharov,
On the rank of the intersection of free subgroups in virtually free groups,
J. Algebra, 418 (2014), 29-43.
\arXiv 1301.3115

\end

\iffalse %!!!!!!!

\Th Mineyev
\rm([Mi12b], part 1.8).
% теоремы 6 в препринте.
Let group $G$, containing
free
% у Минеева отсутствует слово "свободные", но подразумевается видимо
subgroup $A$ and $B$,
acts on some ordered graph $X$,
preserving order and freely on edges (i.e. stabiliser each edge
is trivial). Let $X$ contains
subforest $\cal A$ and $\cal B$, invariant with respect to
subgroups $A$ and $B$,
respectively; and action group $A$ on $\cal A$ and
group $B$ on $\cal B$ freely and cocompact (i.e. has only finite
number orbits vertices and edges).
Then
$
\sum\limits_{u\in S} \r(A^u\cap B)\le\r(A)\cdot\r(B),
$
where summation extends on some
(any)
system $S$
representatives double cosets $AgB$ group $G$ in subgroups~$A$
and $B$.

% Тут явно что-то не то. Должен быть множитель в неравенстве.
% Или "лес" надо заменить на "дерево"..

\fi %!!!!!!!!